\documentclass[11pt]{amsart}
\usepackage{amsmath}
\usepackage{amssymb}
\usepackage{amsthm}
\usepackage{sty1}
\def\bloneg{\mathrm{L}^1(G)}
\def\bltwog{\mathrm{L}^2(G)}
\def\fal#1{\mathrm{A}(#1)}
\def\falg{\mathrm{A}(G)}
\def\falbg#1{\mathrm{A}_{#1}(G)}
\def\fsal#1{\mathrm{B}(#1)}
\def\fsalo#1{\mathrm{B}_0(#1)}
\def\fsalg{\mathrm{B}(G)}
\def\genlin#1{\mathrm{GL}(#1)}
\def\measg{\mathrm{M}(G)}
\def\spine#1{\mathrm{A}^*(#1)}
\def\spinebg#1{\mathrm{A}^*_{#1}(G)}
\def\spineg{\mathrm{A}^*(G)}
\def\tnq#1{\mathcal{T}_{nq}(#1)}
\def\tnqg{\mathcal{T}_{nq}(G)}
\def\vn#1{\mathrm{VN}(#1)}
\def\vnsg{\mathrm{VN}^*(G)}

\begin{document}
%
%
\theoremstyle{plain}
\newtheorem{theorem}{Theorem}[section]
\newtheorem{lemma}[theorem]{Lemma}
\newtheorem{corollary}[theorem]{Corollary}
\newtheorem{proposition}[theorem]{Proposition}
\theoremstyle{definition}
\newtheorem{definition}[theorem]{Definition}
\theoremstyle{remark}
\newtheorem{remark}[theorem]{Remark}
\newtheorem{example}[theorem]{Example}
%
%
\title[Operator amenability of $\fsalg$]{Operator Amenability of Fourier--Stieltjes algebras, II}
\author{Volker Runde}
\address{Department of Mathematical and Statistical Sciences, 
University of Alberta, Edmonton, AB, Canada T6G 2G1}
\thanks{Volker Runde's research supported by NSERC under grant no.\ 227043-04.}
\email{vrunde@ualberta.ca}
\author{Nico Spronk}
\address{Department of Pure Mathematics, University of Waterloo,
Waterloo, ON, Canada N2L 3G1}
\thanks{Nico Spronk's research supported by NSERC under grant no.\ 312515-05.}
\email{nspronk@math.uwaterloo.ca}
\begin{abstract}
We give an example of a non-compact, locally compact group $G$ such that its 
Fourier--Stieltjes algebra $\fsalg$ is operator ame\-nab\-le. Furthermore, we 
characterize those $G$ for which $\spine{G}$---the spine of $\fsalg$ as
introduced by M.\ Ilie and the second named author---is operator amenable
and show that $\spine{G}$ is operator weakly amenable for each $G$.
\end{abstract}
\maketitle
\section{Introduction}
\subsection{History and context}
Let $G$ be a locally compact group, and let $\bloneg$ and $\measg$ denote its 
group and measure algebra, respectively. In \cite{eymard}, P.\ Eymard
introduced the Fourier algebra $\falg$ and the Fourier--Stieltjes algebra 
$\fsalg$ of $G$. In the framework of Kac algebras \cite{enockschwartz},
$\falg$ and $\bloneg$ as well as $\fsalg$ and $\measg$ can be viewed as dual
to one another: this duality generalizes the well known dual group
construction for abelian groups.
\par
It is a now classical theorem of B.\ E.\ Johnson \cite{johnsonM} that
$\bloneg$ is amenable if and only if $G$ is an amenable
group. On the other hand, $\falg$ is an amenable Banach algebra
only if $G$ has an abelian subgroup of finite index \cite{forrestr,runde}.  
In order to obtain the appropriate statement dual to Johnson's theorem, 
we first need to recognize that $\bloneg$---as the predual of the von Neumann
algebra $\mathrm{L}^\infty(G)$---is canonically equipped with an operator 
space structure. In \cite{ruan}, Z.-J.\ Ruan modified Johnson's notion of 
Banach algebraic amenability from \cite{johnsonM} by considering only 
completely bounded module actions and derivations and 
obtained the notion of operator amenability. Since $\mathrm{L}^\infty(G)$ is 
abelian, the canonical operator space structure of $\bloneg$ is 
$\max \bloneg$, so all bounded maps from $\bloneg$ 
are automatically completely bounded; consequently, $\bloneg$ is operator
amenable if and only if it is amenable. As the predual of the group
von Neumann algebra $\mathrm{VN}(G)$, the Fourier algebra $\falg$ also
carries a natural operator space structure, and in \cite{ruan}, Ruan showed
that $\falg$ is operator amenable if and only if $G$ is amenable.  
We further note that Johnson \cite{johnson91} proved that 
$\bloneg$ is always weakly amenable (see \cite{despicg} for a simpler proof)
whereas $\falg$ is known to be weakly amenable only if $G$ if the component 
of the identity is abelian \cite{forrestr}. On the other hand, $\falg$ is 
always operator weakly amenable \cite{spronk,samei1,samei2}.
\par
The questions for which $G$ the measure algebra $\measg$ is amenable or
weakly amenable, respectively, were eventually settled by H.\ G.\ Dales, 
F.\ Ghahramani and A.\ Ya.\ Helemskii \cite{dalesgh}: $\measg$ is
amenable if and only if $G$ is discrete and amenable, and it is weakly 
amenable if and only if $G$ is discrete. In abelian group duality---and
more generally in Kac algebra theory---, the property dual to discreteness
is compactness. Thus, parallel to the $\bloneg$-$\falg$ situation, one 
is led to expect that $\fsalg$ is operator (weakly) amenable if and only if 
$G$ is compact, as was conjectured in both \cite{rundes}
and \cite{spronk} (see also \cite[Problem 9]{helemskii}). 
In fact, it was shown in \cite{rundes} that
$\fsalg$ is operator amenable \emph{with operator amenability constant
$C<5$} if and only if $G$ is compact. At the time \cite{rundes} was written,
the authors felt that the condition imposed on the amenability constant 
was unnecessary.
\par
In the present article, we refute our previous conjecture by exhibiting
examples of non-compact groups $G$ for which $\fsalg$ is operator amenable.
Moreover, we show that the operator amenability constant of $\fsalg$ for
those $G$ is precisely $5$: this shows that the estimate for this constant
from \cite{rundes} cannot be improved. All of our examples, which are 
taken from work by L.\ Bagget \cite{bagget} and G.\ Mauceri and M.\ A.\
Picardello \cite{maucerip}, are separable Fell groups with countable dual
spaces.
\par
The groups $G$ for which $\fsalg$ is amenable were characterized in 
\cite{forrestr}: they are precisely those compact $G$ with an abelian 
subgroup of finite index. In analogy with the corresponding result for 
$\measg$, one might conjecture that $\fsalg$ is weakly amenable exactly
when $G$ is compact with an abelian connected component of the identity.
Our examples show that that this natural conjecture is false too.
\par
Furthermore, we establish some amenability
results for the spine $\spineg$ of $\fsalg$, which was introduced and studied 
by M.\ Ilie and the second-named author in 
\cite{ilies1}.
\subsection*{Acknowledgment}
The second author is pleased to acknowledge a stimulating
discussion he had with Zhiguo Hu, about the role played by
Fell groups in generating counterexamples to natural conjectures about
$\fsalg$.
\subsection{Definitions and notation}
The Fourier and Fourier--Stieltjes algebra, $\falg$ and $\fsalg$, of a 
locally compact group $G$ were introduced in \cite{eymard}.  
the Fourier--Stieltjes algebra $\fsalg$ is the space of all
coefficient functions of weakly operator continuous unitary representations
on Hilbert spaces,
i.e., $\fsalg=\{G \ni s\mapsto\inprod{\pi(s)\xi}{\eta}:(\pi,\fH)\in\Sigma_G, \,
\xi,\eta\in\fH\}$, where $\Sigma_G$ denotes the class of such
representations. Using direct sums and tensor products of representations,
it can be verified that $\fsalg$ is an algebra of functions.
Moreover, $\fsalg$ is the dual space of the
enveloping $\mathrm{C}^\ast$-algebra of $G$, and under this norm is, in fact, 
a Banach algebra. The left regular representation $(\lam,\bltwog)$
is defined through left translation on the space of square integrable
functions with respect to left Haar measure. The Fourier algebra $\falg$ is 
the space of all coefficients of $\lam$, and is a closed ideal in $\fsalg$.
Moreover, $\fsalg$ is the predual of the enveloping von Neumann algebra,
$\mathrm{W}^*(G)$, and $\falg$ is the predual of the group
von Neumann algebra, $\vn{G} := \lambda(G)''$.
\par
In its capacity as the predual of $\mathrm{W}^*(G)$, $\fsalg$
is an \emph{operator space} by \cite{blecher}. Our standard reference
for operator spaces and completely bounded maps is \cite{effrosrB}.
Taking the adjoint of the multiplication map $m_0:\fsalg\otimes\fsalg\to
\fsalg$, one obtains a $^\ast$-homomorphism---and thus a complete 
contraction---$m_0^*:\mathrm{W}^*(G) \to
\mathrm{W}^*(G)\wbar{\otimes}\mathrm{W}^*(G)$, where $\wbar{\otimes}$ denotes 
the von Neumann tensor product. Consequently, $m_0$ extends to a complete 
contraction $m:\fsalg\what{\otimes}\fsalg\to\fsalg$, where $\what{\otimes}$
stands for the operator space projective tensor product
\cite{blecherp,effrosr}. All in all, $\fsalg$ a \emph{completely
contractive Banach algebra}. It is clear that every closed subalgebra of 
$\fsalg$---such as $\falg$---is also a completely contractive Banach algebra.
Note that that $\fsalg\what{\otimes}\fsalg$ is canonically completely 
isometrically isomorphic to a closed subspace of $\fsal{G\cross G}$, and that 
$\falg\what{\otimes}\falg\cong\fal{G\cross G}$ holds completely isometrically 
isomorphically \cite{effrosr}.
\par
If $\fA$ is a completely contractive Banach algebra, a completely
bounded $\fA$-bimodule is an operator space $\fV$, which is a module
for which the module maps $\fA \times \fV \ni (a,v)\mapsto a\mult v$ and 
$\fV \times \fA \ni (v,a)\mapsto v\mult a$
extend to completely bounded maps $\fA\what{\otimes}\fV\to\fV$ and
$\fV\what{\otimes}\fA\to\fA$. Dual modules of completely bounded
$\fA$-bimodules with the dual space operator space structure and the dual 
action are also completely bounded $\fA$-bimodules. A completely
contractive Banach algebra $\fA$ is said to be \emph{operator
amenable} if, for every completely bounded $\fA$-bimodule $\fV$, 
every completely bounded derivation $D:\fA\to\fV^*$ is inner.
\par
The concept of a bounded approximate diagonal can be readily
adapted from the Banach algebra context \cite{johnson72} to the operator
space setting: a completely contractive Banach algebra is amenable 
if, and only if, it admits a \emph{completely bounded approximate diagonal},
i.e., a bounded net $(d_\alp)_\alp$ in $\fA\what{\otimes}\fA$ for which
\[
a\mult d_\alp-d_\alp\mult a\overset{\alp}{\longrightarrow}0
\quad\aand\quad m(d_\alp)a\overset{\alp}{\longrightarrow}a
\]
for each $a\in\fA$, where $a\mult(b\otimes c)=(ab)\otimes c$ and
$(b\otimes c)\mult a=b\otimes(ca)$ for $b, c \in \fA$,
and $m:\fA\what{\otimes}\fA\to\fA$ is the multiplication map.
We say that $\fA$ has \emph{operator amenability constant $C$}, or is
\emph{$C$-operator amenable}, if $C$ is the largest number for which
$\limsup_\alp\norm{d_\alp}\geq C$ for any completely bounded approximate
diagonal. (This notion is again adapted from the corresponding notion for 
Banach algebras, which was developed in \cite{johnson94} specifically
to address non-amenability of Fourier algebras for certain compact groups.)
\par
Following \cite{badecd}, we say that a (completely contractive) commutative
Banach algebra $\fA$ is \emph{(operator) weakly amenable} if
every (completely) bounded
derivation $D:\fA\to\fV$ into a (completely) bounded
symmetric bimodule---i.e., satisfying $a\mult v=v\mult a$ for $a\in\fA$ and
$v\in\fV$---is zero.
\section{Some operator amenable Fourier-Stieltjes algebras}
\label{sec:fsalamen}
For any prime number $p$, let $\Que_p$ denote the field of $p$-adic numbers, 
which is a locally compact field. It is defined to be the completion of the 
rational numbers $\Que$ by the $p$-adic valuation $|\cdot|_p$, which is a 
multiplicative, non-Archimedian valuation, i.e., satisfies 
$|rs|_p=|r|_p|s|_p$ and $|r+s|_p\leq\max\{|r|_p,|s|_p\}$ for $r,s \in \Que_p$. 
 We then let $\Oh_p :=\{r\in\Que_p:|r|_p\leq 1\}$,
the \emph{$p$-adic integers}, which is a compact open subring of $\Que_p$.  
The multiplicative group of $\Oh_p$ is $\Tee_p :=\{r\in\Que_p:|r|_p=1\}$. 
The family of sets $\{p^k\Oh_p\}_{p=0}^\infty$ forms a basis of neighborhoods 
of $0$, so that $\Que_p$ is totally disconnected.
\par
Let $\genlin{n,\Oh_p}$ denote the multiplicative group of
$n\cross n$ matrices with entries in $\Oh_p$ and determinant of
valuation $1$.  This compact group acts on the vector space $\Que_p^n$
by matrix multiplication, and we set
\[
  G_{p,n}:=\genlin{n,\Oh_p}\ltimes\Que_p^n.
\]
For $n=1$, this is the group $\Tee_p\ltimes\Que_p$ of \cite{bagget}.
It is reasonable to call $G_{p,n}$ the ``$n$th rigid $p$-adic motion group''.
\par
In \cite{maucerip}, it was shown---using the ``Mackey machine''---that
the dual space of $\what{G}_{p,n}$---the set of all (equivalence classes) of 
irreducible continuous unitary representations of $G_{p,n}$---is countable.
In fact each group $G_{p,n}$ is of the form $G=K\ltimes A$ where
\begin{enumerate}
\item $K$ is a compact group acting on an abelian group $A$, with
each of the groups separable, and
\item  the dual space $\what{G}$ is countable and decomposes as
$\what{K}\comp q\sqcup\{\lam_k\}_{k=1}^\infty$ where $\what{K}$ is the
discrete dual space of $K$, $q:G\to K$ is the quotient map,
and each $\lam_k$ is a subrepresentation of the left regular representation.
\end{enumerate}
\par
The following was proven for $G_{p,1}$ independently in \cite{mauceri} and 
\cite{walter}. A proof for general $G_{p,n}$ can be obtained in a similar
way. For the reader's convenience, we give a proof for general groups of
the form $G = K \ltimes A$ satisfying (1) and (2).
\begin{proposition} \label{KAprop}
For $G=K\ltimes A$ as above, $\fsalg=\fal{K}\comp q\oplus_{\ell^1}
\falg$ holds.
\end{proposition}
\begin{proof}
Let $u\in\fsalg$, so that $u=\inprod{\pi(\cdot)\xi}{\eta}$ for some
$(\pi,\fH) \in \Sigma_G$ and $\xi,\eta\in\fH$. 
By \cite[Theorem 4.5]{taylor} and (2) above, $\pi$ is totally decomposable.
We may thus write we write
\[
 \pi=\bigoplus_{\sig\in\what{K}} \alp_\sig\mult\sig\comp q\oplus
 \bigoplus_{k=1}^\infty \beta_k\mult\lam_k,
\]
where $\alp_\sig$ and $\beta_k$ are multiplicity constants. For $\sig \in
\what{K}$ and $k \in \mathbb{N}$, let $P_\sig$ and $P_k$ denote the orthogonal 
projection associated with $\alp_\sig\mult\sig\comp q$ and 
$\beta_k\mult\lam_k$, respectively. We then obtain for $s\in G$
that
\[
  u(s)=\sum_{\sig\in\what{K}}
  \inprod{\alp_\sig\mult\sig\comp q(s)P_\sig\xi}{P_\sig\eta}
  +\sum_{k=1}^\infty\inprod{\beta_k\mult\lam_k(s)P_k\xi}{P_k\eta}.
\]
By \cite[(3.13) Corollaire]{arsac}, which uses standard von Neumann
algebra techniques, this is an $\ell^1$-direct sum, i.e.,
\[
  \norm{u}=\sum_{\sig\in\what{K}}
  \norm{\inprod{\alp_\sig\mult\sig\comp q(\cdot)P_\sig\xi}{P_\sig\eta}}
  +\sum_{k=1}^\infty\norm{\inprod{\beta_k\mult\lam_k(\cdot)P_k\xi}{P_k\eta}}.
\]
Each $\inprod{\alp_\sig\mult\sig\comp q(\cdot)P_\sig\xi}{P_\sig\eta}$ lies
in $\fsal{K}\comp q$, and each 
$\inprod{\beta_k\mult\lam_k(\cdot)P_k\xi}{P_k\eta}$ belongs to $\falg$. 
All in all, we see that $u=u_1+u_2$ with $\| u \| = \| u_1 \| + \| u _2 \|$,
where $u_1\in\fal{K}\comp q$ and $u_2\in\falg$.
\end{proof}
\par
Since $\genlin{n,\Oh_p}$ and $G_{p,n}$ are both disconnected, it 
follows from \cite[Section 5]{forrest}, that 
$\fal{\genlin{n,\Oh_p}} \comp q \cong \fal{\genlin{n,\Oh_p}}$
and $\fal{G_{p,n}}$ are both generated by idempotents.  Hence, it
follows from Proposition \ref{KAprop} that $\fsal{G_{p,n}}$ is generated by 
idempotents as well.
\par
From the remarks following the proof of \cite[Theorem 1.4]{badecd}, we
thus obtain:
\begin{corollary} \label{wacor}
For each $n \in \mathbb{N}$, the Fourier--Stieltjes algebra
$\fsal{G_{p,n}}$ is weakly ame\-nab\-le.
\end{corollary}
\par
Since weak amenability implies operator weak amenability, Corollary
\ref{wacor} already shows that---contrary to what one might expect in 
view of \cite{dalesgh}---there are non-compact, locally compact groups with
an operator weakly ame\-nab\-le Fourier--Stieltjes algebra.
\par
We shall now see that $\fsal{G_{p,n}}$ is even operator amenable. In fact,
we will work again in the slightly more general setting of Proposition
\ref{KAprop}.
\begin{theorem} \label{theo:fsalgopamen}
For $G=K\ltimes A$ as in Proposition \emph{\ref{KAprop}}, 
$\fsalg=\fal{K}\comp q\oplus_{\ell^1}\falg$ is 
operator amenable with operator amenability constant $5$.
\end{theorem}
\begin{proof}
This proof is adapted from \cite[Theorem 3.1(i)]{johnson96}.
\par
Since $K$ and $G=K\ltimes A$ are amenable groups, it follows from
(an inspection of) \cite{ruan} that $\fal{K}\comp q\cong\fal{K}$ and $\falg$ 
are each $1$-operator amenable.  
Let
$(u_\alp)_{\alp\in A}$ be a norm $1$ completely
bounded approximate diagonal for $\falg$ and let $(v_\beta)_{\beta\in B}$ be 
such for $\fal{K}$.  
Since $K$ is compact, we can arrange for $v_\beta(k,k)=1$ for all $k\in K$.
Again, $m:\fsalg\what{\otimes}\fsalg\to\fsalg$ denotes the completely 
contractive multiplication map.
\par
Let $\Gamma=A\cross A^A\cross B$ be the product directed set 
\cite[p.\ 69]{kelley}, let $e_\alp=m(u_\alp)$ for each $\alp \in A$, and,
for each $\gam=\bigl(\alp,(\alp'_{\alp''})_{\alp''\in A}, \beta\bigr)$
in $\Gamma$, set
\[
  \begin{split}
  w_\gam &= \bigl((1-e_\alp)\otimes(1-e_{\alp'_\alp})
             +u_{\alp'_\alp}\bigr)v_\beta\comp(q\cross q) \\
         &= \bigl(1\otimes 1-e_\alp\otimes 1-1\otimes e_{\alp'_\alp}
             +e_\alp\otimes e_{\alp'_\alp}+u_{\alp'_\alp}\bigr)
             v_\beta\comp(q\cross q).
  \end{split}
\]
Since
\begin{multline*}
\fsalg\what{\otimes}\fsalg = \fal{K}\comp q\what{\otimes}\fal{K}\comp q\, \\
\oplus_{\ell^1}\,
\fal{K}\comp q\what{\otimes}\falg 
\oplus_{\ell^1}\,\falg\what{\otimes}\fal{K}\comp q\,\oplus_{\ell^1}\,
\falg\what{\otimes}\falg
\end{multline*}
by \cite[Lemma 3.1]{rundes}, we see that
\[
\norm{w_\gam} \leq \left(1+\norm{e_\alp}+\norm{e_{\alp'_\alp}}
+\norm{e_\alp\otimes e_{\alp'_\alp}} +\norm{u_{\alp'_\alp}}
\right)\norm{v_\beta}=5
\]
for each $\gam\in \Gamma$.  
\par
We shall now verify that $(w_\gam)_{\gam\in\Gamma}$ 
is a completely bounded approximate diagonal for $\fsalg$.
\par
We first check that $(w_\gam)_{\gam\in\Gamma}$ is asymptotically central
for the right and left module actions. For $u\in\falg$, we have 
\[
  \begin{split} 
  \lefteqn{\norm{u\mult w_\gam-w_\gam\mult u}} & \\
  & \quad \leq\norm{(u-ue_\alp)\otimes(1-e_{\alp'_\alp})
  -(1-e_\alp)\otimes(u-ue_{\alp'_\alp})+u\mult u_{\alp'_\alp}-u_{\alp'_\alp}
  \mult u} \\
  & \quad \leq 2\norm{u-ue_\alp}+2\norm{u-ue_{\alp'_\alp}}+
  \norm{u\mult u_{\alp'_\alp}-u_{\alp'_\alp}\mult u},
  \end{split}
\]
which can be made arbitrarily small for sufficiently large $\alp$ and 
$\alp'_\alp$. For $v\in\fal{K}$, we have 
\[
\norm{(v\comp q)\mult w_\gam-w_\gam\mult (v\comp q)}
\leq 5\norm{v\mult v_\beta-v_\beta w}
\]
which can be made arbitrarily small for sufficiently large choices of $\beta$.
For general $w=v\comp q+u \in \fsalg$ with $v \in \fal{K}$ and $u \in \falg$, 
it is thus clear that $\norm{w\mult w_\gam- w_\gam\mult w}$ can be made 
arbitrarily small for sufficiently large $\gam\in\Gamma$.
\par
Next, we check that $\bigl(m(w_\gam)\bigr)_{\gam\in\Gamma}$ is an 
approximate identity for $\fsalg$. Note that, for $\gam=
\bigl(\alp,(\alp'_{\alp''})_{\alp''\in A},\beta\bigr) \in \Gamma$ we have
\[
  m(w_\gam)=1+(e_{\alp'_\alp}-1)e_\alp.
\]
For $v\in\fal{K}\comp q$, we then obtain
\[
  \lim_\gam m(w_\gam)v-v=\lim_{\alp}\lim_{\alp'}(e_{\alp'}-1)e_\alp v=0 
\]
because $\lim_{\alp'}(e_{\alp'}-1)e_\alp v=0$ for each $\alp$ as $e_\alp 
v\in\falg$. A similar calculation shows that $\lim_\gam m(w_\gam)u=u$ for
$u\in\falg$. 
\par
Consequently, $( w_\gamma )_{\gamma \in \Gamma}$ is indeed a completely
bounded approximate diagonal for $\fsalg$, so that the operator
amenability constant $C$ of $\fsalg$ can at most be $5$. Since $G$ is not
compact, $C<5$ cannot occur by \cite[Theorem 3.2]{rundes}. Hence, $C = 5$
must hold.
\end{proof}
\begin{remark}
The operator amenability of $\fsalg$ in Theorem \ref{theo:fsalgopamen} can
easily be obtained by observing that $\falg$ is an operator amenable ideal in 
$\fsalg$ with operator amenable quotient
$\fsalg/\falg \cong \fal{K}$ and then applying the operator space
analog of \cite[Theorem 2.3.7]{rundeB}. The disadvantage in doing this,
however, is that it yields no information on the operator amenability constant of
$\fsalg$.
\end{remark}
\begin{remark}
In \cite{rundes}, we introduced, for an arbitrary locally compact group $G$, 
a decomposition $\fsalg = \mathrm{A}_{\mathcal{F}}(G) \oplus_{\ell^1}
\mathrm{A}_{\mathcal{PIF}}(G)$, which can be interpreted as dual to the
decomposition of $\measg$ into the discrete and the continuous measures.
For non-discrete $G$, it is well known that there are continuous measures
in $\measg \setminus \bloneg$, and we conjectured that, at least for
amenable $G$, the inclusion $\falg \subset \mathrm{A}_{\mathcal{PIF}}(G)$
is proper unless $G$ is compact
\cite[p.\ 681, Remarks (2)]{rundes}. Theorem \ref{theo:fsalgopamen}
shows that this conjecture is false.
\end{remark}
\section{Operator amenability of the spine}\label{sec:spine}
\subsection{The spine of $\fsalg$}
The main theorem of the previous section is actually a particular case of a 
more general result. The groups $G_{p,n}$ are all non-compact
amenable groups for which $\spine{G_{p,n}}=\fsal{G_{p,n}}$ holds, where
$\spine{G_{p,n}}$ is the spine of $\fsal{G_{p,n}}$ as defined in 
\cite{ilies1}. In this section, we study amenability properties for spines.
\par
We recall the definition of the \emph{spine} of a Fourier--Stieltjes 
algebra $\fsalg$ of a locally compact group $G$ below. 
Full details are presented 
in the article \cite{ilies1}.
\par
Let $\tnqg$ denote the family of all group topologies $\tau$ on $G$ with 
the following properties:
\begin{itemize}
\item The completion $G_\tau$ of $G$ with respect to the left 
uniformity generated by $\tau$ is a locally compact group. (This completion
is unique up to homeomorphic isomorphism, and it is the same completion as
gained from the right uniformity.)
\item If $\tau_{ap}$ is the coarsest topology making the almost periodic
compactification map $\eta:G\to G^{ap}$ continuous, then 
$\tau\supseteq\tau_{ap}$.
\end{itemize}
We call such topologies \emph{non-quotient locally precompact 
topologies}. The family $\tnqg$ is a semilattice, i.e., a commutative,
idempotent semigroup, under the operation $(\tau_1,\tau_2)\mapsto
\tau_1\vee\tau_2$, where $\tau_1\vee\tau_2$ is the coarsest topology
which is simultaneously finer than both $\tau_1$ and $\tau_2$.
In particular, this semilattice is unital with unit $\tau_{ap}$.
\par
For each $\tau\in\tnqg$ we let $\eta_\tau:G\to G_\tau$ be the natural map
into the completion. Then the Fourier algebra $\fal{G_\tau}$ is completely
isometrically isomorphic to the subalgebra $\falbg{\tau}=
\fal{G_\tau}\comp\eta_\tau$ of $\fsalg$.  We have that
$\falbg{\tau_1}\cap\falbg{\tau_2}=\{0\}$ if $\tau_1\not=\tau_2$
and $\norm{u_{\tau_1}+u_{\tau_2}}=\norm{u_{\tau_1}}+\norm{u_{\tau_2}}$ for $u_{\tau_j}\in
\falbg{\tau_j}$ ($j=1,2$), in this case.  The spine is then the algebra
\[
\spineg=\text{$\ell^1$-}\!\!\!\!\bigoplus_{\tau\in\tnqg}\falbg{\tau},
\]
which is graded over $\tnqg$, i.e., $u_{\tau_1}u_{\tau_2}\in\falbg{\tau_1 \vee\tau_2}$
for $u_{\tau_j}\in \falbg{\tau_j}$ ($j=1,2$).
\par
As in \cite[Section 6.2]{ilies1}, we can calculate that $\tnq{G_{p,n}}=
\{\tau_{ap},\tau_{p,n}\}$ for any of the groups
$G_{p,n}$, where $\tau_{p,n}$ is the given
topology on $G_{p,n}$.  Moreover, the quotient map
$q:G_{p,n}\to\genlin{n,\Oh_p}$ is the almost periodic compactification
map. Consequently, $\fsal{G_{p,n}}=\spine{G_{p,n}}$ holds.
\subsection{Operator amenability of $\spineg$}
For our discussion of the operator amenability of $\spineg$, we introduce some
auxiliary notation. For any $F \subseteq \tnqg$, let $\langle F \rangle$ denote
the sublattice of $\tnqg$ it generates, and define
\[
  \spinebg{F}=\text{$\ell^1$-}\!\bigoplus_{\tau\in\langle F\rangle}
  \falbg{\tau}.
\]
Note that, if $F\subseteq \tnqg$ is finite, then so is $\langle F\rangle$.
\par
Since $\tnqg$ is finite for each of the groups $G_{p,n}$, 
the following lemma extends (the qualitative part of) 
Theorem \ref{theo:fsalgopamen}.
\begin{proposition}\label{spineprop}
Let $G$ be an amenable, locally compact group, and let $F\subseteq \tnqg$ be
finite. Then $\spinebg{F}$ is operator amenable.
\end{proposition}
\begin{proof}
We shall prove that $\spinebg{F}$ is amenable by using induction on 
$|F|$. 
\par
Suppose that $|F| = 1$, so that $F=\{\tau \}$ for some $\tau \in \tnqg$.
Since $G$ is amenable, so is $G_\tau$ by \cite[Proposition 1.2.1]{rundeB},
which implies that $\spinebg{F} \cong \fal{G_\tau}$ is operator amenable
by \cite{ruan}.
\par
Now suppose that $|F|>1$. Fix $\tau \in F$, let $F' := F \setminus
\{ \tau \}$, and set
\[
  I_{\tau,F} := \mathrm{A}^\ast_{\tau\vee F'}(G) + \mathrm{A}_\tau(G),
\]
where $\tau\vee F'=\{\tau\vee\tau':\tau'\in F'\}$. Then $I_{F,\tau}$
is an ideal in $\spinebg{F}$ containing $\mathrm{A}^\ast_{\tau\vee F'}(G)$
as an ideal. Since $\mathrm{A}^\ast_{\tau\vee F'}(G)$ is operator amenable by
the induction hypothesis, and since $I_{\tau,F} / 
\mathrm{A}^\ast_{\tau\vee F'}(G)$ is either $\mathrm{A}_\tau(G)$ or $\{ 0 \}$, 
we conclude from the completely bounded analogue of 
\cite[Theorem 2.3.10]{rundeB} that $I_{\tau,F}$ is operator amenable.
Since $\spinebg{F} = \spinebg{F'} + I_{\tau,F}$ and $\spinebg{F'}$ is operator
amenable by induction hypothesis, a similar argument yields
the operator amenability of $\spinebg{F}$.
\end{proof}
\par
It is immediate from Proposition \ref{spineprop} that $\spineg$ is
operator amenable if $G$ is amenable and $\tnqg$ is finite. We shall see
that these are the only conditions under which $\spineg$ can be operator
amenable.
\par
For the the following lemma, note that, by linearity and continuity,
the product of $\tnqg$ extends to $\ell^1(\tnqg)$ turning it into a
Banach algebra. Since the canonical operator space structure of
$\ell^1(\tnqg)$ is $\max \ell^1(\tnqg)$, this Banach algebra is canonically
completely contractive.
\begin{lemma} \label{spinelem}
Let $G$ be a locally compact group. Then the map
\[
  \Pi : \spineg \to \ell^1(\tnqg), \quad ( u_\tau)_{\tau \in \tnqg} \mapsto
  \sum_{\tau \in \tnqg} u_\tau(e) \, \delta_\tau
\]
is a complete quotient map and an algebra homomorphism.
\end{lemma}
\begin{proof}
Note that $\vnsg := \text{$\ell^\infty$-} \bigoplus_{\tau\in\tnqg}\vn{G_\tau}$ 
is the dual space of $\spineg$. For each $\tau \in \tnqg$, let $p_\tau
\in \vnsg$ be the central projection corresponding to the identity element
of $\vn{G_\tau}$. Then
\[
  \ell^\infty(\tnqg) \to \vnsg, \quad (\lambda_\tau)_{\tau \in \tnqg}
  \mapsto ( \lambda_\tau \, p_\tau )_{\tau \in \tnqg}
\]
is a normal $^\ast$-monomorphism and the adjoint of $\Pi$. This
show that $\Pi$ is indeed a complete quotient map.
\par
To see that $\Pi$ is multiplicative, let $\tau_1, \tau_2 \in \tnqg$,
and let $u_{\tau_j} \in \mathrm{A}_{\tau_j}(G)$ for $j = 1,2$. It follows
that $u_{\tau_1} u_{\tau_2} \in \mathrm{A}_{\tau_1 \vee \tau_2}(G)$ and thus
\[
  \Pi(u_{\tau_1} u_{\tau_2}) = 
  u_{\tau_1}(e) u_{\tau_2}(e) \delta_{\tau_1 \vee \tau_2} =
  u_{\tau_1}(e) \, \delta_{\tau_1} \; u_{\tau_2}(e) \, \delta_{\tau_2} =
  \Pi(u_{\tau_1}) \, \Pi( u_{\tau_2}). 
\]
By linearity and continuity, this proves the multiplicativity of $\Pi$.
\end{proof}
\begin{theorem} \label{spinethm}
Let $G$ be a locally compact group. Then $\spineg$ is operator
amenable if and only if $G$ is amenable and $\tnqg$ is finite.
\end{theorem}
\begin{proof}
The ``if'' part is provided by Proposition \ref{spineprop}.       
\par
For the ``only if'' part, suppose that 
$\spineg$ is operator amenable. The central projection $p_G \in \vnsg$ 
corresponding to the identity operator
in $\vn{G}$ forms a completely contractive projection onto $\falg$
via the predual action, $\spineg \ni u\mapsto p_G\mult u$.  Thus
$\falg$ is a completely complemented ideal in $\spineg$
and hence is operator amenable 
by the completely bounded analog of \cite[Theorem 2.3.7]{rundeB}.
Therefore, by \cite{ruan}, $G$ is amenable.
\par
Since $\spineg$ is operator amenable, so is its quotient $\ell^1(\tnqg)$,
and since the canonical operator space structure of $\ell^1(\tnqg)$ is
$\max \ell^1(\tnqg)$, it follows that $\ell^1(\tnqg)$ is amenable in the
purely Banach algebraic sense. From \cite[Theorem 2.7]{gronbaek}, we
conclude that $\tnqg$ is finite.
\end{proof}
\begin{example}
Using computations from \cite[Section 6]{ilies1}, we obtain 
that $\spineg$ is operator amenable for $G$ being any one of
the following groups: the real numbers $\Ree$, the integers $\Zee$, the 
Euclidean motion groups $\mathrm{M}(n)=\mathrm{SO}(n)\ltimes\Ree^n$ for
$n \in \mathbb{N}$, the $ax+b$ group, or $\Que_p$, where $p$ is any prime.  
On the other hand, the spine fails to be operator amenable for
any of the groups $\Ree^n$ or $\Zee^n$ with $n\geq 2$, 
for $\Que$ as a discrete group, and for any non-amenable group.
\end{example}
\begin{remark}
For any locally compact group $G$, let $\fsalo{G}$ denote the closed ideal of 
$\fsalg$ consisting of functions vanishing at $\infty$. 
For the Euclidean motion groups, 
it is known (see, for example, the discussion on \cite[p.\ 10]{chou})
that $\fsal{\mathrm{M}(n)}=\fal{\mathrm{SO}(n)}\comp q\oplus_{\ell^1}
\fsalo{\mathrm{M}(n)}$. We suspect that
$\fsalo{G}$ is never operator amenable when it is properly larger than
$\falg$ (see \cite{figatalamanca} for situations in which this is known to
be the case). This would entail that $\fsal{\mathrm{M}(n)}$ cannot be operator 
amenable.
\end{remark}
\begin{remark}
In view of Theorems \ref{theo:fsalgopamen} and \ref{spinethm}, we a prepared
to make the conjecture that, for a locally compact group $G$, the 
Fourier--Stieltjes algebra $\fsalg$ is operator amenable if and only if
$\fsalg = \spineg$, $G$ is amenable, and $\tnqg$ is finite.  
\end{remark}
\par
We note:
\begin{corollary}
Let $G$ be a locally compact group. Then $\spineg$ is amenable
if and only if $G$ has an abelian subgroup group of finite index
and $\tnqg$ is finite.
\end{corollary}
\begin{proof}
If $\spineg$ is amenable, it is operator amenable, so that $\tnqg$ 
must be finite. Since $\falg$ is a complemented ideal of $\spineg$, it
must be amenable, too. Hence, $G$ has an abelian subgroup of finite index
by \cite{forrestr,runde}.
\end{proof}
\subsection{Operator weak amenability of $\spineg$}
In contrast to Theorem \ref{spinethm}, we have the following:
\begin{proposition}
Let $G$ be a locally compact group $G$. Then $\spineg$ is operator
weakly amenable.  
\end{proposition}

\begin{proof}
Let $\fV$ be an completely bounded symmetric $\spineg$-bimodule, and 
let $D:\spineg\to\fV$ be a completely bounded derivation. Then
\[
  D=\sum_{\tau\in\tnqg}D|_{\falbg{\tau}}.
\]
holds. Since $\falbg{\tau}\cong\fal{G_\tau}$ is operator weakly amenable
by \cite{spronk}, it follows that $D|_{\falbg{\tau}}=0$ for each
$\tau \in \tnqg$ and thus
$D = 0$.
\end{proof}
\begin{remark}
It is not clear at all for which locally compact groups $G$, the spine
$\spineg$ might be weakly amenable (in the original Banach algebraic sense).
The spine \emph{is} weakly amenable for any compact group with an abelian 
connected component of the identity by \cite[Theorem 3.3]{forrestr} and also 
for any of the groups $G_{p,n}$ by Corollary \ref{wacor}. However,
there are compact groups $K$ for which $\fal{K}$ is not weakly amenable 
\cite{johnson94}. If we let $G$ be any discrete 
group for which $G^{ap}$ admits such a group $K$ as a quotient, then
$\spineg$ appears not to be weakly amenable.  Indeed, $\falbg{\tau_{ap}}\cong
\fal{G^{ap}}$ is a quotient of $\spineg$.  Furthermore 
$\fal{G^{ap}}$ contains an isometric copy of $\fal{K}$. Thus, it appears
unlikely that $\fal{G^{ap}}$ is weakly amenable, and we conjecture the
same for $\spineg$. 
\end{remark}
\bibliographystyle{amsplain}
\bibliography{fsalgopamen}

\providecommand{\bysame}{\leavevmode\hbox to3em{\hrulefill}\thinspace}
\providecommand{\MR}{\relax\ifhmode\unskip\space\fi MR }
\providecommand{\MRhref}[2]{%
  \href{http://www.ams.org/mathscinet-getitem?mr=#1}{#2}
}
\providecommand{\href}[2]{#2}
\begin{thebibliography}{10}

\bibitem{arsac}
G.~Arsac, \emph{Sur l'espace de {B}anach engendr\'{e} par les coefficients
  d'une repr\'{e}sentation unitaire}, Pub. D\'{e}p. Math. Lyon \textbf{13}
  (1976), no.~2, 1--101.

\bibitem{badecd}
W.~G. Bade, P.~C. Curtis, and H.~G. Dales, \emph{Amenability and weak
  amenability for {B}eurling and {L}ipschitz algebras}, Proc. London Math. Soc.
  \textbf{55} (1987), no.~3, 359--377.

\bibitem{bagget}
L.~Bagget, \emph{A separable group having discrete dual space is compact}, J.
  Funct. Anal. \textbf{10} (1972), 131--148.

\bibitem{blecher}
D.~P. Blecher, \emph{The standard dual of an operator space}, Pacific Math. J.
  \textbf{153} (1992), 15--30.

\bibitem{blecherp}
D.~P. Blecher and V.~I. Paulsen, \emph{Tensor products of operator spaces}, J.
  Funct. Anal. \textbf{99} (1991), 262--292.

\bibitem{chou}
C.~Chou, \emph{Minimally weakly almost periodic groups}, J. Funct. Anal.
  \textbf{36} (1980), 1--17.

\bibitem{dalesgh}
H.~G. Dales, F.~Ghahramani, and A.~Ya. Helemskii, \emph{The amenability of
  measure algebras}, J. London Math. Soc. \textbf{66} (2002), 213--225.

\bibitem{despicg}
M.~Despi\'{c} and F.~Ghahramani, \emph{Weak amenablility of group algebras of
  locally compact groups}, Canad. Math. Bull. \textbf{37} (1994), no.~2,
  165--167.

\bibitem{effrosr}
E.~G. Effros and Z.-J. Ruan, \emph{On approximation properties for operator
  spaces}, International J. Math. \textbf{1} (1990), 163--187.

\bibitem{effrosrB}
\bysame, \emph{Operator spaces}, London Math. Soc., New Series, vol.~23,
  Claredon Press, Oxford Univ. Press, New York, 2000.

\bibitem{enockschwartz}
M.~Enock and J.-M. Schwartz, \emph{Kac algebras and duality of locally compact
  groups}, Springer Verlag, Berlin--Heidelberg--New York, 1992.

\bibitem{eymard}
P.~Eymard, \emph{L'alg\`{e}bre de {F}ourier d'un groupe localement compact},
  Bull. Soc. Math. France \textbf{92} (1964), 181--236.

\bibitem{figatalamanca}
A.~Fig\`{a}-Talamanca, \emph{Positive definite functions which vanish at
  infinity}, Pacific. J. Math. \textbf{69} (1977), 355--363.

\bibitem{forrest}
B.~E. Forrest, \emph{Fourier analysis on coset spaces}, Rocky Mountain J. Math.
  \textbf{28} (1998), 173--189.

\bibitem{forrestr}
B.~E. Forrest and V.~Runde, \emph{Amenability and weak amenability of the
  {F}ourier algebra}, Math. Z. \textbf{250} (2004), 731--744.

\bibitem{gronbaek}
N.~Gr{\o}nb{\ae}k, \emph{Amenability of discrete convolution algebras, the
  commutative case}, Pacific J. Math. \textbf{143} (1990), 243--249.

\bibitem{helemskii}
A.~Ya. Helemskii, \emph{Some aspects of topological homology since 1995: a
  survey}, Banach algebras and their applications (A.~T.-M. Lau and V.~Runde,
  eds.), Contemporary Mathematics, vol. 363, American Mathematical Society,
  2004, pp.~145--179.

\bibitem{ilies1}
M.~Ilie and N.~Spronk, \emph{The spine of a {F}ourier--{S}tieltjes algebra},
  Proc. London Math. Soc. (accepted); see Ar{X}iv {\tt math.FA/0405063}, 2005.

\bibitem{johnson72}
B.~E. Johnson, \emph{Approximate diagonals and cohomology of certain
  annihilator {B}anach algebras}, Amer. J. Math. \textbf{94} (1972), 685--698.

\bibitem{johnsonM}
\bysame, \emph{Cohomology in banach algebras}, Mem. Amer. Math. Soc., vol. 127,
  American Mathematical Society, 1972.

\bibitem{johnson91}
\bysame, \emph{Weak amenability of group algebras}, Bull. London Math. Soc.
  \textbf{23} (1991), 281--284.

\bibitem{johnson94}
\bysame, \emph{Non-amenability of the {F}ourier algebra of a compact group}, J.
  London Math. Soc. \textbf{50} (1994), 361--374.

\bibitem{johnson96}
\bysame, \emph{Symmetric amenability and the nonexistence of {L}ie and {J}ordan
  derivations}, Math. Proc. Cambridge Phil. Soc. \textbf{120} (1996), 455--473.

\bibitem{kelley}
J.~L. Kelley, \emph{General topology}, University Series in Higher Mathematics,
  Van Nostrand, New York, 1955.

\bibitem{mauceri}
G.~Mauceri, \emph{Squre integrable representations and the {F}ourier algebra of
  a unimodular group}, Pacific J. Math. \textbf{73} (1977), 143--154.

\bibitem{maucerip}
G.~Mauceri and M.~A. Picardello, \emph{Non-compact unimodular groups with
  purely atomic {P}lancherel measures}, Proc. Amer. Math. Soc. \textbf{78}
  (1980), 77--84.

\bibitem{ruan}
Z.-J. Ruan, \emph{The operator amenability of ${A(G)}$}, Amer. J. Math.
  \textbf{117} (1995), 1449--1474.

\bibitem{rundeB}
V.~Runde, \emph{Lectures on amenability}, Lecture Notes in Mathematics, vol.
  1774, Springer Verlag, Berlin--Heidelberg--New York, 2002.

\bibitem{runde}
\bysame, \emph{The amenability constant of the {F}ourier algebra}, Proc. Amer.
  Math. Soc. \textbf{134} (2006), 1473--1481.

\bibitem{rundes}
V.~Runde and N.~Spronk, \emph{Operator amenability of {F}ourier-stieltjes
  algebras}, Math. Proc. Camb. Phil. Soc. \textbf{136} (2004), 675--686.

\bibitem{samei1}
E.~Samei, \emph{Approximately local derivations}, J. London Math. Soc.
  \textbf{71} (2005), 759--778.

\bibitem{samei2}
\bysame, \emph{Hyper-{T}auberian algebras and weak amenability of
  {F}ig\`{a}-{T}alamanca-{H}erz algebras}, J. Funct. Anal. \textbf{231} (2006),
  195--220.

\bibitem{spronk}
N.~Spronk, \emph{Operator weak amenability of the {F}ourier algebra}, Proc.
  Amer. Math. Soc. \textbf{130} (2002), 3609--3617.

\bibitem{taylor}
K.~F. Taylor, \emph{Geometry of {F}ourier algebras and locally compact groups
  with atomic unitary representations}, Math. Ann. \textbf{262} (1983),
  183--190.

\bibitem{walter}
M.~E. Walter, \emph{On a theorem of {F}ig\`{a}-{T}alamanca}, Proc. Amer. Math.
  Soc. \textbf{60} (1976), 72--74.

\end{thebibliography}
\end{document}